\documentclass[12pt]{amsart}
\usepackage{amsmath}
\usepackage{amssymb}
\usepackage{amsthm}
\usepackage{mathrsfs}
\usepackage{fullpage}
\usepackage{color}

\newcommand{\Aut}{\operatorname{Aut}}
\newcommand{\cue}{\mathbb{Q}}
\newcommand{\fpbar}{\overline{\mathbb{F}}_p}
\newcommand{\fptwo}{\mathbb{F}_{p^2}}
\newcommand{\GL}{\operatorname{GL}}
\newcommand{\GSp}{\operatorname{GSp}}
\newcommand{\GU}{\operatorname{GU}}
\newcommand{\zed}{\mathbb{Z}}

\newtheorem{thm}{Theorem}
\newtheorem{cor}[thm]{Corollary}
\theoremstyle{definition}
\newtheorem*{rem}{Remarks}
\newtheorem*{ack}{Acknowledgements}

\begin{document}
\title{Upper bound on the number of systems of Hecke eigenvalues for
  Siegel modular forms (mod $p$)}
\author[A. Ghitza]{Alexandru Ghitza}
\address{
  Department of Mathematics and Statistics\\
  McGill University\\
  805 Sherbrooke W\\
  Montr\'eal, Qu\'ebec\\
  H3A 2K6, CANADA
}
\email{aghitza@alum.mit.edu}

\begin{abstract}
  We derive an explicit upper bound for the number
  $\mathscr{N}(g,N,p)$ of systems of Hecke eigenvalues coming from
  Siegel modular forms (mod $p$) of dimension $g$ and level $N$
  relatively prime to $p$.  In the special case of elliptic modular
  forms ($g=1$), our result agrees with recent work of G. Herrick.
\end{abstract}

\maketitle

\section{Introduction}
We prove the following
\begin{thm}\label{thm:bound}
  Fix an integer $g\geq 1$, a prime $p$ and an integer $N\geq 3$ not
  divisible by $p$.  The number of systems of Hecke eigenvalues coming
  from Siegel modular forms (mod $p$) of dimension $g$ and level $N$
  satisfies
  \begin{equation*}
    \mathscr{N}(g,N,p)\leq c_g\cdot\#\GSp_{2g}(\zed/N\zed)\cdot
    p^{(g+2)(g-1)/2}\cdot(p^2-1)\cdot\prod_{j=1}^g
    \left(p^j+(-1)^j\right),
  \end{equation*}
  where
  \begin{equation*}
    c_g=(-1)^{\frac{g(g+1)}{2}}2^{-g}\prod_{j=1}^g
    \zeta(1-2j)=\frac{1}{2^{2g}g!}\prod_{j=1}^g B_{2j}.
  \end{equation*}
\end{thm}

\begin{cor}[Asymptotics]
  \ 
  \begin{enumerate}
  \item Fix $g$ and $p$, then
    \begin{equation*}
      \mathscr{N}(g,N,p)=O\left(N^{2g^2+g+1}\right).
    \end{equation*}
  \item Fix $g$ and $N$, then
    \begin{equation*}
      \mathscr{N}(g,N,p)=O\left(p^{g^2+g+1}\right).
    \end{equation*}
  \end{enumerate}
  The constants are effectively computable.
\end{cor}
\begin{proof}
  Part (a) follows from the fact that the algebraic group $\GSp_{2g}$
  has dimension $2g^2+g+1$.
  
  
  Part (b) is obvious.
\end{proof}

Combined with Theorem 1.1 of \cite{Ghitza2004a}, Theorem~\ref{thm:bound}
gives
\begin{cor}[Algebraic modular forms]
  Let $B/\cue$ be the quaternion algebra ramified at $p$ and $\infty$.
  The number of systems of Hecke eigenvalues coming from algebraic
  modular forms (mod $p$) of level $N$ on $\GU_g(B)$ satisfies the
  inequality of Theorem~\ref{thm:bound}.
\end{cor}

\begin{rem}
  \begin{enumerate}
  \item In the case of classical modular forms ($g=1$) with
    $\Gamma(N)$ level structure we get $\mathscr{N}(1,N,p)=O(p^3)$;
    the author recently learned of an upper bound with the same
    asymptotics obtained by Graham Herrick for $\Gamma_0(N)$ and
    $\Gamma_1(N)$ level structures.  The previously known bound was
    $O(p^4)$ and can be found in~\cite{Jochnowitz1982a}.
  \item A more general approach in the context of algebraic modular
    forms is presented in Section 7 of~\cite{Gross1996}.
  \end{enumerate}
\end{rem}

\begin{ack}
  The author is indebted to A. J. de Jong and M. Lieblich for helpful
  conversations and to B. Gross for his inspiring work on algebraic
  modular forms.  This research was funded by the Centre
  Interuniversitaire en Calcul Math\'ematique Alg\'ebrique (CICMA) and
  the Fonds Qu\'ebecois de la Recherche sur la Nature et les
  Technologies (FQRNT).
\end{ack}

\section{The setup}
We fix the dimension $g$, the prime $p$ and the level $N$ not
divisible by $p$.  We define (see Sections 2.2 and 2.3
of~\cite{Ghitza2004a} for details)
\begin{eqnarray*}
  \Sigma&=&\{[A,\lambda]:A\text{ superspecial abelian variety over
  }\fpbar, \lambda\text{ principal polarization}\},\\
  \Sigma(N)&=&\{[A,\lambda,\alpha]:\text{ as above}, \alpha\text{
    symplectic level $N$ structure on }(A,\lambda)\},\\
  \tilde{\Sigma}(N)&=&\{[A,\lambda,\alpha,\eta]:\text{ as above},
  \eta\text{ hermitian basis of invariant differentials over
  }\fptwo\text{ on }(A,\lambda)\}.
\end{eqnarray*}
The notation $[x]$ stands for the isomorphism class of the object $x$.

It follows from the proof of Theorem 4.5 of~\cite{Ghitza2004a} that the
restriction of Siegel modular forms to the superspecial locus
$\Sigma(N)$ induces a bijection on the sets of systems of Hecke
eigenvalues.  Therefore the number that we want to estimate is the
number of systems of Hecke eigenvalues occurring in the spaces of
superspecial modular forms
\begin{equation*}
  S_\tau(N)=\left\{f:\tilde{\Sigma}(N)\to
    W_\tau\,|\,f([A,\lambda,\alpha,M\eta])=
    \tau(M)^{-1}f([A,\lambda,\alpha,\eta])\text{
      for all }M\in\GU_g(\fptwo)\right\},
\end{equation*}
as the weight $\tau:\GU_g(\fptwo)\to\GL(W_\tau)$ runs through the set
of irreducible representations of $\GU_g(\fptwo)$ over $\fpbar$.

Note that for a fixed class $[A,\lambda,\alpha]\in\Sigma(N)$, a
superspecial modular form $f\in S_\tau(N)$ is completely determined by
the value it takes on $[A,\lambda,\alpha,\eta]$ for any choice of
$\eta$, and therefore
\begin{equation*}
  \dim S_\tau(N)\leq \#\Sigma(N)\cdot\dim W_\tau.
\end{equation*}
We get the following upper bound for the number we are interested in:
\begin{equation*}
  \mathscr{N}(g,N,p)\leq \#\Sigma(N)\cdot\sum_{\tau}\dim W_\tau.
\end{equation*}

\section{The cardinality of $\Sigma(N)$}
We write
\begin{equation}
  \#\Sigma(N)=\sum_{[A,\lambda,\alpha]\in\Sigma(N)}
  1=\sum_{[A,\lambda]\in\Sigma}
  \#\{[A,\lambda,\alpha]\}.\label{eqn:sigma-n-1}
\end{equation}
So we fix $(A,\lambda)$ and we want to count the number of level $N$
structures $\alpha$, up to $\fpbar$-isomorphism.  If we ignore the
isomorphisms, there are precisely $\#\GSp_{2g}(\zed/N\zed)$ level $N$
structures.  But by definition, two level $N$ structures $\alpha$ and
$\alpha'$ on $(A,\lambda)$ are isomorphic if and only if
$\alpha'=\alpha\circ\varphi$ for some automorphism $\varphi$ of
$(A,\lambda)$.  Hence we can continue our computation from
(\ref{eqn:sigma-n-1}) as follows:
\begin{eqnarray}\label{eqn:sigma-n}
  \#\Sigma(N)&=&\sum_{[A,\lambda]\in\Sigma}
  \#\{[A,\lambda,\alpha]\}=\sum_{[A,\lambda]\in\Sigma}
  \frac{\#\GSp_{2g}(\zed/N\zed)}{\#\Aut(A,\lambda)}\notag \\
  &=&c_g\cdot\#\GSp_{2g}(\zed/N\zed)\cdot
  \prod_{j=1}^g\left(p^j+(-1)^j\right).
\end{eqnarray}
Here we used the following mass formula due to Ekedahl (page 159
of~\cite{Ekedahl1987}, see also Corollary 9.5 in~\cite{Geer1999}):
\begin{equation*}
  \sum_{[A,\lambda]\in\Sigma}
  \frac{1}{\#\Aut(A,\lambda)}=c_g\cdot\prod_{j=1}^g\left(p^j+(-1)^j\right).
\end{equation*}

\section{Representations of $\GU_g(\fptwo)$}
We need to estimate the sum of the dimensions of the irreducible
representations of $\GU_g(\fptwo)$ defined over $\fpbar$.  For this we
use the theory of finite groups with split $(B,N)$-pairs, as explained
in~\cite{Curtis1970}.


The cardinality of the group $\GU_g(\fptwo)$ is 
\begin{equation*}
  \#\GU_g(\fptwo)=p^{g(g-1)/2}(p-1)\prod_{j=1}^g\left(p^j-(-1)^j\right).
\end{equation*}
Therefore the size of a $p$-Sylow subgroup is $p^{g(g-1)/2}$, and so
Corollary 5.11 of~\cite{Curtis1970} says that
\begin{equation*}
  \dim W_\tau\leq p^{g(g-1)/2}\quad\quad\text{for any irreducible }\tau.
\end{equation*}
The rank of $\GU_g(\fptwo)$ is $g$, and so by 
Proposition 6.1 of~\cite{Gross1996}
we know that the number of irreducible representations of
$\GU_g(\fptwo)$ over $\fpbar$ is $p^{g-1}(p^2-1)$.  Therefore we get
the inequality
\begin{equation}\label{eqn:dim}
  \sum_\tau \dim W_\tau\leq
  p^{g-1}(p^2-1)p^{g(g-1)/2}=p^{(g+2)(g-1)/2}\cdot(p^2-1).
\end{equation}

\section{The end}
It remains to put (\ref{eqn:sigma-n}) and (\ref{eqn:dim}) together to
get
\begin{equation*}
  \mathscr{N}(g,N,p)\leq c_g\cdot\#\GSp_{2g}(\zed/N\zed)\cdot
  p^{(g+2)(g-1)/2}\cdot(p^2-1)\cdot\prod_{j=1}^g
  \left(p^j+(-1)^j\right),
\end{equation*}
which is precisely the content of Theorem~\ref{thm:bound}.

\end{document}